\documentclass[10pt, a4paper]{article}
\begin{document}

\title {Extension of the Erberlein-\u{S}mulian Theorem to Normed
        Spaces}
\author{ Wha Suck Lee\\
         Department of Mathematics, University Of Pretoria, Pretoria, South Africa\\
       }
\date {November 18. 2002}
\maketitle

\begin{center}
\small{wlee@postino.up.ac.az}
\end{center}

\begin{abstract}

   The Erberlein-\u{S}mulian theorem asserts that for {\it {complete}} normed
   spaces, that is, Banach spaces, a subset is weak compact if and only if it is weak
   sequentially compact. In this paper it is shown that the completeness of the normed
   space is not necessary for the above mentioned result.

\end{abstract}

\bigskip
\setcounter {footnote} {0}
\section {Introduction.}

   Since a normed space is trivially a metric space,
   a subset is norm topology compact if and only if it is
   norm topology sequentially compact. As far as the concepts of
   closed convex subsets, bounded subsets or continuous linear functionals are
   concerned, they are the same both in the norm and the weak topology.
   Just as in the case of the norm topology where norm
   topology compactness is equivalent to norm topology sequential
   compactness, it is to be expected that weak sequential compactness
   should be equivalent to weak compactness in {\it {all}} normed
   spaces.\\

   \noindent Here we establish this result by showing that the
   steps in the proof of the Erberlein- \u{S}mulian theorem, found in J. Diestel
   [1], Chapter 3, pp 18 - 20, hold for all normed spaces, without the assumption of
   completeness.

\pagebreak

\section {Weak compactness.}

Here we establish the weak compactness of a bounded subset in a
normed space, not necessarily complete, by using Alaoglu's theorem: \\

{\bf Alaoglu's theorem.}{ \it {For any normed linear space X,
$B_{X^ {*}}$, the closed unit ball in the dual $X^{*}$, is
$weak^{*}$ compact. Consequently, $weak^{*}$ closed and (norm)
bounded subsets of the normed space $X^{*}$ are $weak^{*}$
compact. }} \\

Let A be a bounded subset of a normed space X. To show $\bar{A} ^
{weak}$ is $weak$ compact, that is, A is relatively $weak$
compact, we need to see how $\bar{A} ^ {weak}$ looks like.

Consider A as J[A] where J is the canonical embedding of X into
$X^{**}$. Let us look at $\bar{A} ^ {weak^{*}}$ in the bidual
$X^{**}$. $\bar{A} ^ {weak}$ is identical to $\bar{A} ^
{weak^{*}}$, provided that the elements in $X^{**}$ found in
$\bar{A} ^ {weak^{*}}$ are precisely point evaluational
functionals, that is, in X. Furthermore, the $weak^{*}$ and the
$weak$ topologies coincide. Hence, once we show that $\bar{A} ^
{weak^{*}}$ is $weak^{*}$ compact by using Alaoglu's theorem, then
$\bar{A} ^ {weak}$ will also be weak compact.

\bigskip

{\bf {Proposition 1.}} {\it{ $\bar{A} ^ {weak^{*}}$ is $weak^{*}$
 compact.}}\\

{\bf Proof.} By Alaoglu's theorem, to show $\bar{A} ^ {weak^{*}}$
is $weak^{*}$ compact it suffices to show that $\bar{A} ^
{weak^{*}} $  is (norm) bounded in the normed space $X^{**}$. By
Lemma 1 (please see the appendix), applied to the normed linear
space $X^{**}$, it is equivalent to show that $\bar{A} ^
{weak^{*}}$ is weakly bounded (as defined in Lemma 1) in the
bidual $X^{**}$. We therefore show that:\\

($\forall f $ $\in X^{***}$) ($\sup_{x^{**} \in \bar {A} ^
{weak^{*} } } $ $\{ \mid f (x^{**}) \mid$ $\leq K_{f} \}$) where
$K_{f}$ is a constant which
depends on $f$.  \\

Consider an arbitrary f $\in X^{***}$. Then $\parallel f \parallel
\ \leq \ R_{f} $ for some $R_{f}  \geq 0$; that is, f resides in a
ball in $X^{***}$ of radius $R_{f}$. By Goldstine's theorem
(please see the appendix), $R_{f} B_{X^{*}}$ is $weak^{*}$ dense
in $R_{f} B_{X^{***}}.$ Consequently, there exists a net of point
evaluational functionals in  $X^{***},$ indexed by $\delta$, $\
\Theta_ {x^{*}_{\delta}} \in R_{f} B_{X^{*}}$, such that:

\begin{equation}
 \Theta_ {x^{*}_{\delta}} \stackrel
{weak^{*}}{\longrightarrow}f
\end{equation}

 \bigskip

 \noindent Equivalently,

 \begin{equation}
  (\forall x^{**} \in X^ {**} ) \ x^{**} (x^{*}_ \delta )\longrightarrow f(x^{**})
 \end{equation}

\bigskip
\noindent
 where $\parallel  \Theta_ {x^{*}_{\delta}} \parallel
 \   = \ \parallel {x^{*}_{\delta}} \parallel \  \leq \ R_{f}.
$
\bigskip

In particular, taking $x^{**}_{0} \in \bar{A} ^ {weak^{*}}$, we
have the convergence of the following complex valued net, indexed
by $\delta$, $ \ x^{**}_{0} (x^{*}_\delta )$ :

 \begin{equation}
    \fbox {
    $x^{**}_{0} (x^{*}_\delta ) \longrightarrow f(x^{**}_{0}) \ where \
     \parallel {x^{*}_{\delta}} \parallel \  \leq  \   R_{f}.$   }
 \end{equation}

\pagebreak
\noindent Since $x^{**}_{0} \in \bar {A} ^ {weak^{*}}$ where we
treat A as J[A], that is,  A consists entirely of point
evaluational  functionals in $X^{**}$, there exists a net of point
evaluational functionals in $A \subset X^{**},$ indexed by
$\gamma$, $\ \Theta_ {x_{\gamma}}$, such that:

\begin{equation}
\Theta_ {x_{\gamma}} \stackrel {weak^{*}}{\longrightarrow}
x^{**}_{0} \  where \parallel \Theta_{x_{\gamma}} \parallel =
\parallel {x_{\gamma}}\parallel \ \leq \ r.
\end{equation}

\noindent where r is the  radius of the ball in the bidual
$X^{**}$ in which A resides in: A is bounded, by hypothesis.\\

In (4), the $weak^{*}$ convergence is in $X^{**}$, where the
$weak^{*}$ dual of $X^{***}$ is $X^{*}$. Hence, (4) is equivalent
to:


\begin {center}
($\forall \Theta_{x^{*}}) \ (x^{*} \in X^{*}$) \  $\Theta_{x^{*}}
(\Theta _ {x_{\gamma}}) \longrightarrow \Theta_{x^{*}}
(x^{**}_{0}))$
\end{center}
\noindent or
\begin{equation}
 (\forall x^{*} \in X^{*})\ \   x^{*} (x_{\gamma})
\longrightarrow x^{**}_{0} (x^{*} )
\end{equation}
\noindent \\ In particular, applying (5) to each $x^{*}_{\delta}$
of
formula (3): $x^{*}_{\delta} \in X^{*}$ , we have: \\
\begin{equation} 
 \fbox {
 $\forall x^{*}_{\delta}$ of formula (3) \ \ \  $x^{*}_{\delta} (x_{\gamma})
\longrightarrow x^{**}_{0} (x^{*}_{\delta} ) $ \ \ \ (index
$\gamma$ runs, the index $\delta$ fixed )}
 \end{equation}  
\noindent where $\parallel x_{\gamma} \parallel \leq r$ and
$\parallel x^{*}_{\delta}
\parallel \leq R_{f}$.\\

\noindent Therefore by (6) and (3) respectively,
\begin{center}
 $(\forall \delta) \  x^{*}_{\delta} (x_{\gamma}) \longrightarrow
 \ \
 x^{**}_{0} (x^{*}_{\delta})$ ($\gamma$ runs only)
\end{center}
\noindent and
\begin{center}
 $x^{**}_{0} (x^{*}_{\delta}) \ \ \longrightarrow \ \  f(x^{**}_{0})$ ($\delta$ runs)
\end{center}

\noindent But\\

\begin{center}
$\mid x^{*}_{\delta} (x_{\gamma}) \mid \ \leq \ \parallel
x^{*}_{\delta}\parallel \parallel x_{\gamma}\parallel \ \leq \
R_{f}
* r \ \ \  \forall \delta$
\end{center}

\bigskip
Therefore the positive valued numerical net $\mid x^{*}_{\delta}
(x_{\gamma}) \mid$ is bounded by $R_{f} * r$. By the continuity of
the norm,$\mid \cdot \mid$, in the complex field, $\mid
x^{*}_{\delta} (x_{\gamma}) \mid \longrightarrow
 \ \
 \mid x^{**}_{0} (x^{*}_{\delta}) \mid $. Since $\mid x^{**}_{0} (x^{*}_{\delta})
 \mid$ is the limit of the positive valued numerical net $\mid x^{*}_{\delta}
(x_{\gamma}) \mid$ which resides in the closed unit ball of radius
$R_{f} * r $, the limit $\mid x^{**}_{0} (x^{*}_{\delta})
 \mid$  also resides in this ball. Repeating this argument on
 formula (3), $\mid f(x^{**} _ {0}) \mid$, being the limit of the
 net $\mid x^{**}_{0} (x^{*}_{\delta}) \mid$ which resides in the
 closed ball of radius $R_{f} * r$, will also reside here too.
 Therefore, $\mid f(x^{**}_{0}) \mid \leq R_{f} * r $.

\begin{flushright} Q.E.D \end{flushright}

\pagebreak

{\bf Proposition 2.} { \it {Let us look at $\bar{A} ^ {weak^{*}}$
in the bidual $X^{**}$. $\bar{A} ^ {weak}$ is identical to
$\bar{A} ^ {weak^{*}}$, provided that the elements in $X^{**}$
found in $\bar{A} ^ {weak^{*}}$ are precisely point evaluational
functionals, that is, in X. Furthermore, the $weak^{*}$ and the
$weak$ topologies coincide.}}

\bigskip

Let us note that if $X^{**}$ has a locally convex topology, the
topology relative to any subspace, in particular X, is also
locally convex and the fundamental neighbourhood filter of the
subspace is the trace of the fundamental neighbourhood system
filter of $X^{**}$ on X.

Let $\mathcal{U}_{0}$ denote the fundamental neighbourhood system
of the $weak^{*}$ topology of $X^{**}$. The typical neighbourhoods
are sets of the form $\{ x^{**} \in X^{**} \mid \\ \mid
x^{**}(x^{*}_{i}) \mid \leq \epsilon , \ \  i = 1, \ldots ,n  \}$.

\noindent Then \\ \begin{eqnarray*} X \cap \{ x^{**} \in X^{**}
\mid \ \ \mid x^{**} (x^{*}_{i}) \mid \leq \epsilon &&,i = 1,
\ldots ,n \} \\ & = &  \{ \Theta_{x} \in X^{**}\mid \ \mid
\Theta_{x}(x^{*}_{i}) \mid \leq \epsilon, i = 1, \ldots ,n \}
\\   &  = & \{ \Theta_{x} \in X^{**} \mid \ \mid (x^{*}_{i})(x)
\mid \leq \epsilon, \ \ i = 1, \ldots ,n \}
\end{eqnarray*}

\noindent which are precisely the typical neighbourhoods of X
endowed with the weak topology.

Since the elements in $X^{**}$ found in $\bar{A} ^ {weak^{*}}$ are
precisely point evaluational functionals of $X^{**}$, that is, in
X, $\bar{A}^{weak^{*}} \cap X = \bar{A}^{weak^{*}}$. Since the
$weak^{*}$ topology of $X^{**}$ relativised with respect to X is
the weak topology, $\bar {A}^{X} = \bar {A} ^{weak}$. By Lemma 3
(please see the appendix), the conclusion follows.

\begin{flushright} Q.E.D \end{flushright}

\section {Necessary condition in normed spaces.}

Here we prove the necessary condition in the Eberlein-\u{S}mulian
Theorem in normed spaces, that is, without using completeness.
Namely, we show that a relatively weak compact
subset of a {\it{normed}} space is relatively weak sequentially compact.\\

Let A be a relatively weak compact subset of the normed space
X.\\Let ($a_{n}$) be a sequence taken from A. Let [$a_{n}$] denote
the closed linear span of ($a_{n}$). Since the norm closure of a
subspace coincides with the weak closure, we conclude that
[$a_{n}$] is weak closed in X. Therefore, $A \cap [ a_{n} ] $ is
relatively weak compact (compactness is closed hereditary) in the
separable normed space [$a_{n}$].

Since the dual of a separable normed space [$a_{n}$] has a
countable total set by Lemma 6, letting $A \cap [ a_{n} ] $ and
[$a_{n}$] play the roles of A and X, respectively, in Lemma 5, we
infer that $\overline {A \cap [ a_{n} ]} ^{weak} $ is metrizable
in the weak topology of [$a _{n}$]. In metric spaces, compactness
and sequential compactness coincide. Therefore, $\overline {A \cap
[ a_{n} ]} ^{weak} $ is sequentially compact in the weak topology
of [$a_{n}$].

\pagebreak

Hence for sequence ($a_{n}) \subset \overline {A \cap [ a_{n} ]}
^{weak} \subset \bar{A}^{weak}$, there exists a subsequence
($a_{n_{k}})$ which converges to a point, $a$, in $\overline {A
\cap [ a_{n} ]} ^{weak} \subset \bar{A}^{weak}$ in the subspace
topology of $[a_{n}]$, where $a \in [a_{n}] $. Since $a$ is a
limit point of  ($a_{n_{k}})$ with respect to the $[a_{n}]$ -
subspace topology, it will also be a limit point of ($a_{n_{k}})$
with respect to the weak topology of X; this follows
from the definition of the subspace topology.\\

\section {Sufficient condition in normed spaces.}

Here we prove the sufficient condition in the Eberlein-\u{S}mulian
Theorem in normed spaces, that is, without using completeness.
Namely, we show that a relatively weak sequentially compact
subset of a {\it{normed}} space is relatively weak compact.\\

Let A be a relatively weak sequentially compact space of X. In
anticipation of the criteria established in section 2, to show
that $\bar{A} ^ {weak}$ is $weak$ compact (i.e A is relatively
$weak$
compact), we need to show that A is norm bounded first.

It suffices to show that $\bar{A} ^ {weak} \supset A$ is weak
bounded since the norm bounded subsets of a normed space X are
precisely the weak bounded subsets of X. This follows by Lemma 9
since $(X, X^{*})$ is a dual pair and both the norm and the
$\sigma (X, X^{*})$ are topologies of this dual pair. Hence by
Lemma 9, the norm and weak bounded subsets are identical.\\

{\bf {Proposition 1.}} {\it{ If $\bar{A}^{weak}$ is a weak
sequentially compact set of X, then $\bar{A}^{weak}$ is weak
bounded or equivalently norm bounded.}}\\

{\bf {Proof.}} Suppose $\bar{A}^{weak}$ is not weak bounded. Then
we show that $\bar{A}^{weak}$ is not weak sequentially compact.

If $\bar{A}^{weak}$ is not weak bounded, then there exists a weak
neighbourhood $U$ of the convex space $(X, \sigma (X, X^{*}))$
such that for each $n \in {\bf{N}}$, there exists an $a_{n} \in
\bar{A}^{weak}$ such that $a_{n} \notin nU$. The range of the
sequence $( a_{n} )$ cannot have a convergent subsequence. Any
convergent subsequence is Cauchy and hence bounded by Lemma 8
(please see the appendix) : there exists a $m \in {\bf{N}} $ such
that the subsequence is a subset of $mU$ . But any subsequence of
$a_{n}$ cannot be absorbed by $mU$.
\begin{flushright} Q.E.D \end{flushright}

Having established the fact that $A \subset \bar{A}^{weak}$ is
bounded, by Section 2, Proposition 1,  we conclude that $\bar{A} ^
{weak^{*}}$ is $weak^{*}$ compact. To show $\bar {A} ^{weak}$ is
weak compact, all that is left is to show that $\bar{A} ^
{weak^{*}}$ resides entirely in X. To do this, we construct a
point evaluational functional that is arbitrarily close to each
$x^{**} \in \bar{A} ^ {weak^{*}}$, enabling us to identify
$x^{**}$ with that point evaluational functional.

The point evaluational functional that will do the job, is derived
from a sequence $(a_{n}) \subset A \subset \bar {A} ^{weak}$, each
$a_{n}$ being regarded as a point evaluational functional. Let us
recall that by hypothesis $\bar {A} ^{weak}$ is weak sequentially
compact. Therefore, $a_{n}$ has a convergent subsequence
$(a_{n_{k}})$. The point evaluational functional that will do the
job  is precisely the limit of the convergent subsequence
$(a_{n_{k}})$ .

\pagebreak

\noindent We now construct the sequence $(a_{n})$.

\subsection {The sequence $(a_{n})$.}

We now introduce two new terminologies:\\

{\bf{\noindent Terminology 1: Norm Optimizing Finite Set of a subspace E, $NOFS_{E}$: }}\\

Let us recall from Lemma 7 (please see the appendix) that for each
finite dimensional subspace of the bidual $X^{**}, E$, we can
associate a {\it {finite}} subset $E ^{0}$ of the unit {\it
{sphere}} of the dual $X^{*}, S_{X^{*}}$ with the following
property:

\begin{equation}
\frac{\parallel x^{**} \parallel}{2} \ \leq max \ \{ |
x^{**}(x^{*})| \  where \  x^{*} \in E^{0}  \}
\end{equation}

\noindent for each $x^{**} \in \ finite \ dimensional \ subspace \
E$.   \\

\noindent [NOTE:$S_{X^{*}} = \{ x^{*} \in X^{*} | \parallel x^{*}
\parallel = 1 \}$ ] \\

Since, $\parallel x^{**}\parallel = sup_{\parallel x^{*}
\parallel \leq \ 1} $ $\{ \mid x^{**}(x^{*}) \mid \}$, the finite set
$E^{0} \subset S_{X^{*}}$ associated with the finite dimensional
subspace E of $X^{**}$ will be called a {\it{norm optimizing
finite set of subspace E}} which we abbreviate as $NOFS_{E}$.\\

{\bf{\noindent Terminology 2: the $(x^{*}_{1},1)$-cluster point of $x^{**}$: }}\\

\noindent Consider $x^{**} \in \bar{A} ^ {weak^{*}}$\\
Let $x^{*}_{1}$ be a member of unit sphere $S_{X^{*}}$. In the
$weak^{*}$ topology, the family of seminorms $\{ |\Theta_{x^{*}}|
\ \ | \ x^{*} \in X^{*}\}$ generates the defining family of
seminorms. Now the $weak^{*}$ neighbourhood of $x^{**}$, generated
by the unit ball of the continuous seminorm
$|\Theta_{x^{*}_{1}}|$, $x^{**} + W^{*}(0_{X^{**}},
\Theta_{x^{*}_{1}}, 1) = \{ y^{**} \in X^{**}| \ \mid (y^{**} -
x^{**})(x^{*}_{1}) \mid < \ 1 \}$, contains a member $a_{1}$ of A.
Consequently,\begin{center} $\mid (x^{**} - a_{1}) (x_{1}^{*})
\mid < \ 1 \ \ \ \  a_{1} \in A \subset X \hookrightarrow X^{**}.$
\end{center}

\noindent We call the point evaluational functional of the bidual
$X^{**}$, $a_{1}$, {\it {the $(x^{*}_{1}, 1)$-cluster point of
$x^{**}$.}}

Similarly, we call the point evaluational functional, say, $a_{2}
\in A \subset X \hookrightarrow X^{**}$ where $a_{2} \in x^{**} +
W^{*}(0_{X^{**}}, \Theta_{x^{*}_{1}},\Theta_{x^{*}_{2}}, 1) \cap A
= \{ y^{**} \in X^{**}| \ \mid (y^{**} - x^{**})(x^{*}_{1}) \mid
and \ \mid (y^{**} - x^{**})(x^{*}_{2}) \mid \ < \ 1 \} \cap A $,
{\it{the $(x^{*}_{1}, x^{*}_{2}, 1)$-cluster point of $x^{**}$. }}
Consequently,\begin{center} $\mid (x^{**} - a_{i}) (x_{i}^{*})
\mid < \ 1 \ \ \ \  a_{i} \in A \subset X \hookrightarrow X^{**}
 where \ \ i = 1, 2.$
\end{center}

\pagebreak

{\bf{\noindent The construction of $(a_{n})$.}}\\

First we start with $x^{**} \in \bar{A} ^ {weak^{*}}$. Choose any
$x^{*}_{1} \in S_{X^{*}}$.

We then let $a_{1}$ be the $(x^{*}_{1}, 1)$-cluster point of
$x^{**}.$ Consequently, \begin{eqnarray*} \mid (x^{**} - a_{1})
(x_{1}^{*}) \mid < \ 1. \end{eqnarray*}

Then let $E_{2} = [x^{**}, x^{**} - a_{1}]:  E_{2}$ is the finite
dimensional subspace of $X^{**}$ spanned by $x^{**}, x^{**} -
a_{1}$. We then associate with $E_{2}$ a norm optimizing finite
subset, $NOFS_{E_{2}} = \{x^{*}_{2}, \ldots, x^{*}_{n(2)}\}
\subset S_{X^{*}}$.

We let $a_{2}$ be the $(NOFS_{E_{2}} \cup \{ x^{*}_{1}
\},\frac{1}{2})$-cluster point of $x^{**}.$  Consequently,

\begin{eqnarray*} \mid (x^{**} - a_{2}) (x_{i}^{*}) \mid < \ \frac{1}{2} \ for
\ i = 1, \ldots ,n(2) \\ \forall y^{**} \in E_{2}, \frac{\parallel
y^{**} \parallel}{2} \ \leq max \ \{ | y^{**}(x^{*}_{i})| \  where
\  i = 1, \ldots ,n(2)
\end{eqnarray*}

Then let $E_{3} = [x^{**}, x^{**} - a_{1} , x^{**} - a_{2}]$:
$E_{3}$ is the finite dimensional subspace of $X^{**}$ spanned by
$x^{**}, x^{**} - a_{1} , x^{**} - a_{2}$. We then associate with
$E_{3}$, a norm optimizing finite subset, $NOFS_{E_{3}} =
\{x^{*}_{n(2) + 1}, \ldots, x^{*}_{n(3)}\} \subset S_{X^{*}}$.

We let $a_{3}$ be the $(NOFS_{E_{3}} \cup NOFS_{E_{2}} \cup \{
x^{*}_{1} \},\frac{1}{3})$-cluster point of $x^{**}.$
Consequently,

\begin{eqnarray*}
\mid (x^{**} - a_{3}) (x_{i}^{*}) \mid < \ \frac{1}{3} \ for \ i =
1, \ldots ,n(3)\\\forall y^{**} \in E_{3}, \frac{\parallel y^{**}
\parallel}{2} \ \leq max \ \{ | y^{**}(x^{*}_{i})| \  where \  i =
1, \ldots ,n(3)
\end{eqnarray*}

Continuing this process inductively, we have for the general case,
n, $E_{n} = [x^{**}, x^{**} - a_{1} \ldots  ,x^{**} - a_{n - 1}
]$: $E_{n}$ is the finite dimensional subspace of $X^{**}$ spanned
by  $x^{**}, x^{**} - a_{1} \ldots  ,x^{**} - a_{n - 1}$ . We then
associate with $E_{n}$, a norm optimizing finite subset,
$NOFS_{E_{n}} = \{x^{*}_{n(n - 1) + 1}, \ldots, x^{*}_{n(n)}\}
\subset S_{X^{*}}$.

We let $a_{n}$ be the $(NOFS_{E_{n}} \cup \ldots \cup NOFS_{E_{2}}
\cup \{ x^{*}_{1} \},\frac{1}{n})$-cluster point of $x^{**}.$
Consequently,
\begin{eqnarray*}
\mid (x^{**} - a_{n}) (x_{i}^{*}) \mid < \ \frac{1}{n} \ for \ i =
1, \ldots ,n(n)\\\forall y^{**} \in E_{n}, \frac{\parallel y^{**}
\parallel}{2} \ \leq max \ \{ | y^{**}(x^{*}_{i})| \  where \  i =
1, \ldots ,n(n)\\ \end{eqnarray*}

Hence we have an ascending chain of subspaces in $X^{**}$ : $E_{1}
\subset E_{2} \subset \ldots \subset E_{n} \subset \ldots$. We,
now give a schematic diagram of the construction.

\pagebreak

\noindent {\bf{The construction of $(a_{n})$ schematically:}} \\

\noindent [ $x^{**} \in \bar{A} ^ {weak^{*}}$]

\noindent Start  $x^{**} \in \bar{A} ^ {weak^{*}}$ (treat $x^{**}$ as a new origin).\\

\begin{flushright}$a_{1}$is the $(x^{*}_{1}, 1)$-cluster point of $x^{**}.$\end{flushright}
\begin{flushright}$\mid (x^{**} - a_{1}) (x_{1}^{*}) \mid < \ 1$\end{flushright}

\noindent $E_{2} = [x^{**}, x^{**} - a_{1}] \longleftrightarrow
NOFS_{E_{2}} = \{x^{*}_{2}, \ldots, x^{*}_{n(2)}\} \hookrightarrow
S_{X^{*}}$

\begin{flushright}$a_{2}$is the $(NOFS_{E_{2}} \cup \{ x^{*}_{1} \},\frac{1}{2})$-cluster point of $x^{**}.$\end{flushright}
\begin{flushright} $\mid (x^{**} - a_{2}) (x_{i}^{*}) \mid < \ \frac{1}{2} \ for \ i = 1, \ldots ,n(2)$\end{flushright}

\begin{flushright} {\it{$\forall y^{**} \in E_{2}, \frac{\parallel
y^{**} \parallel}{2} \ \leq max \ \{ | y^{**}(x^{*}_{i})| \  where
\  i = 1, \ldots ,n(2)$} }
\end{flushright}
\begin{flushright} \end{flushright}

\bigskip

\noindent $E_{3} = [x^{**}, x^{**} - a_{1} , x^{**} - a_{2}]
\longleftrightarrow NOFS_{E_{3}} = \{x^{*}_{n(2) + 1}, \ldots,
x^{*}_{n(3)}\} \hookrightarrow S_{X^{*}}$

\begin{flushright}$a_{3}$is the $(NOFS_{E_{3}} \cup NOFS_{E_{2}} \cup \{ x^{*}_{1} \},\frac{1}{3})$-cluster point of $x^{**}.$\end{flushright}
\begin{flushright} $\mid (x^{**} - a_{3}) (x_{i}^{*}) \mid < \ \frac{1}{3} \ for \ i = 1, \ldots ,n(3)$\end{flushright}

\begin{flushright} {\it{$\forall y^{**} \in E_{3}, \frac{\parallel
y^{**} \parallel}{2} \ \leq max \ \{ | y^{**}(x^{*}_{i})| \  where
\  i = 1, \ldots ,n(3)$} }
\end{flushright}
\begin{flushright} \end{flushright}

 \bigskip
\begin{center} \ \vdots \end{center}

 \noindent $E_{n} = [x^{**}, x^{**} - a_{1} \ldots  ,x^{**}
- a_{n - 1} ] \longleftrightarrow NOFS_{E_{n}} = \{x^{*}_{n(n - 1)
+ 1}, \ldots, x^{*}_{n(n)}\} \hookrightarrow S_{X^{*}}$

\begin{flushright}$a_{n}$is the $(NOFS_{E_{n}} \cup \ldots \cup NOFS_{E_{2}} \cup \{ x^{*}_{1} \},\frac{1}{n})$-cluster point of $x^{**}.$\end{flushright}
\begin{flushright} $\mid (x^{**} - a_{n}) (x_{i}^{*}) \mid < \ \frac{1}{n} \ for \ i = 1, \ldots ,n(n)$\end{flushright}

\begin{flushright} {\it{$\forall y^{**} \in E_{n}, \frac{\parallel
y^{**} \parallel}{2} \ \leq max \ \{ | y^{**}(x^{*}_{i})| \  where
\  i = 1, \ldots ,n(n)$} }
\end{flushright}
\begin{flushright} \end{flushright}

 \bigskip
 \begin{center} \ \vdots \end{center}
\pagebreak

Let us recall that to show $\bar {A} ^{weak}$ is weak compact, all
that was left was to show that $\bar{A} ^ {weak^{*}}$ resides
entirely in X. To do this, we construct a point evaluational
functional that is arbitrarily close to each $x^{**} \in \bar{A} ^
{weak^{*}}$, enabling us to identify $x^{**}$ with that point
evaluational functional.

The point evaluational functional that will do the job, is derived
from the sequence $(a_{n}) \subset A \subset \bar {A} ^{weak}$,
just constructed beforehand; each $a_{n}$ being regarded as a
point evaluational functional. Let us recall that by hypothesis
$\bar {A} ^{weak}$ is weak sequentially compact. Therefore,
$a_{n}$ has a convergent subsequence $(a_{n_{k}})$. The point
evaluational functional that will do the job  is precisely the
limit of the convergent subsequence $(a_{n_{k}})$. We denote this
limit as $x \ or \ \Theta_{x}$. Then the sequence $(a_{n})$ is
frequently in every neighbourhood of $x \ or \ \Theta_{x}$.

We now prove this assertion in the following four steps: \\

{\bf{Proof.}}\\ \\ \indent {\bf {Step 1:}}{\it {$x \in [a_{n}]$.
}}
\\ \\ \indent The closed linear span [$a_{n}$] of the constructed
sequence $(a_{n})$ is weakly closed; therefore, $x \in [a_{n}]$.\\

\indent {\bf {Step 2:}}{\it {$\ x^{**} - x \ is \ in \ the \
weak^{*} \ closed \ linear \ span \ of \ ( x^{**}, x^{**} - a_{1},
x^{**} - a_{2}, \ldots ) $ }}. \\

Let $V = weak^{*} \  closure \ of \ span \ (a_{n})$. Then $x \in \
weak \ closure \ of \ span \ (a_{n}) = [a_{n}] \ (by \ step \ 1)\
\subset V$ (since the $weak^{*}$ topology is coarser than the weak
topology of the bidual $X^{**} \hookleftarrow [a_{n}]$).
Consequently, $x^{**} - x \in x^{**} - V \subset W$ where $W$
denotes the $\ weak^{*} -closed \ linear \ span \ of \ ( x^{**},
x^{**} - a_{1}, x^{**} - a_{2}, \ldots ) $.\\ \\

\noindent {\bf{Proof:}} $v \in V $ iff $\exists \ net \ x_{\delta}
\stackrel{weak^{*}}{\longrightarrow} v$. By the continuity of the
map $L_{x^{**}}: X^{**} \rightarrow X^{**} | z \mapsto x^{**} - z$
with respect to the $weak^{*}$ topology, $\
L_{x^{**}}(x_{\delta})\stackrel{weak^{*}}{\longrightarrow}
L_{x^{**}}(v)$. Equivalently,  $\ x^{**} -
x_{\delta}\stackrel{weak^{*}}{\longrightarrow} x^{**} - v$. Each
member of the net $x_{\delta}$ is of the form $\sum^{n}_{i = 1}
\lambda_{i} a_{i}$ where $a_{i} \in (a_{n})$, since they lie in
$span(a_{n})$. Since $ \lambda_{i} a_{i} =  - \lambda_{i} (x^{**}
- a_{i}) + \lambda_{i} x^{**}$, the translated net $x^{**} -
x_{\delta}$ which has typical members of the form $x^{**} -
\sum^{n}_{i = 1} \lambda_{i} a_{i}$ will belong to $span ( x^{**},
x^{**} - a_{1}, x^{**} - a_{2}, \ldots )$. So, setting $v
= x$, the conclusion follows. \\

\indent {\bf {Step 3:}}{\it {$\forall  y^{**}\ in \ the \ weak^{*}
\ closed \ linear \ span \ of \ ( x^{**}, x^{**} - a_{1}, x^{**} -
a_{2}, \ldots )$
 \begin{equation}  \frac{\parallel y^{**}
\parallel}{2} \ \leq \sup_{m} \  | y^{**}(x^{*}_{m}) | \end{equation}

\noindent $ where \ x_{m}^{*} \in \{ x_{1}^{*} \} \cup \bigcup_{m}
NOFS_{E_{m}}  \ \ and \ \  a_{i} \in$ the constructed sequence $\{ a_{n}\}.$ }} \\

By the ascending chain $E_{1} \subset E_{2} \subset \ldots \subset
E_{n} \subset \ldots$, for any $y^{**} \in span (x^{**}, x^{**} -
a_{1}, \ldots ,x^{**} - a_{n}, \ldots), \ y^{**} \in E_{m} \ for \
some \ m.$ Consequently, \begin{center}$\frac{\parallel y^{**}
\parallel}{2} \ \leq max_{i =
1, \ldots ,n(m)} \ \{ | y^{**}(x^{*}_{i})| \  \   \leq
\frac{\parallel y^{**} \parallel}{2} \ \leq \sup_{n} \ \{ |
y^{**}(x^{*}_{n})| $.\end{center} We now show (10) holds for  all
$y^{**} \ in \ the \ weak^{*} \ closed \ linear \ span \ of \ (
x^{**}, x^{**} - a_{1}, x^{**} - a_{2}, \ldots )$ using a
continuity argument:

\pagebreak

Suppose $y^{**}\ is \ in \ the \ weak^{*} \ closed \ linear \ span
\ of \ ( x^{**}, x^{**} - a_{1}, x^{**} - a_{2}, \ldots )$. Then
$\exists \ net \ x_{\delta} \in span \ of \ (x^{**}, x^{**} -
a_{1}, x^{**} - a_{2}, \ldots) \  | \ x_{\delta} \stackrel
{weak^{*}} {\longrightarrow}  y^{**} $ if and only if
$x_{\delta}(x^{*}) \longrightarrow y^{**}(x^{*})$ for all $x^{*}
\in X^{*}$.Consequently $\parallel x_{\delta} \parallel
\longrightarrow \parallel y^{**} \parallel $ and $\sup_{m \in
{\bf{N}}} |x_{\delta}(x^{*}_{m})| \longrightarrow \sup_{m \in
{\bf{N}}} |y^{**}(x^{*}_{m})|$.

 For each $\forall \delta \ \
\frac{\parallel x_{\delta} \parallel}{2} \leq \sup_{m \in {\bf
{N}}} |x_{\delta} (x^{*}_{m}) | $  where $\ x_{m}^{*} \in \{
x_{1}^{*} \} \cup \bigcup_{m} NOFS_{E_{m}}$. Taking limits (limits
preserve order),$\frac{\parallel y^{**} \parallel}{2} \leq \sup_{m
\in {\bf {N}}} |y^{**} (x^{*}_{m}) | $.\\ \\

{\bf{Step 4}} {\it{Apply Step 3, equation (10) to $x^{**} - x$}}.\\

Since $x^{**} - x \in  \ the \ weak^{*} \  closed \ linear \ span
\ of \ ( x^{**}, x^{**} - a_{1}, x^{**} - a_{2}, \ldots )$, we can
apply (10) by step 3. But \begin{eqnarray} |(x^{**} -
x)(x^{*}_{m})| & = & |(x^{**} - a_{k} + a_{k} - x)(x^{*}_{m})| \nonumber \\
& = & |(x^{**} - a_{k})(x^{*}_{m}) + (a_{k} - x)(x^{*}_{m})| \nonumber \\
& \leq & |(x^{**} - a_{k})(x^{*}_{m})| + |(a_{k} - x)(x^{*}_{m})|
\end{eqnarray}

\noindent for each $a_{k} \in $ constructed sequence $(a_{n})$.\\

Recall the construction of $(a_{n})$ (section 4.1). Then setting
$k \geq p$ for the term $a_{k}$, working in $E_{k} \subset E_{p}$,
for any $m \leq n(p)$, we have:\\

\begin{center} $|(x^{**} - a_{k})(x^{*}_{m}) | \leq \frac{1}{p}$
\end{center}

Recall also that $x$ is a weak cluster point of sequence $(a_{n})
\hookrightarrow X$. In fact, subsequence $(a_{n_{k}})
\stackrel{weak}{\longrightarrow} x$. Therefore, ($X^{*}$ is the
weak dual of $X$)

\begin{equation} x^{*}(a_{n_{k}})
\longrightarrow x^{*}(x)\ \ \  for \ each \ x^{*} \in X^{*}.
\end{equation}

In equation (10), we took $x^{**}, a_{k}, x \in X^{**}$ where
$a_{k}, x$ are the point evaluational functionals
$\Theta_{a_{k}} , \Theta_{x}$. So:\\

\begin{center} $|(a_{k} - x)(x^{*}_{m})| = |x^{*}_{m}(a_{k}) - x^{*}_{m}(x)|$
\end{center}

\noindent can be made as small as possible by taking $k$  high
enough. (by equation (12).) Therefore, $x^{**} - x = 0$ so that
$x^{**} = x \in X$.

\begin{flushright} Q.E.D \end{flushright}
\pagebreak
\setcounter {footnote} {0}
\begin{center} APPENDIX \end{center}

 {\bf Lemma 1.}{ \it {For any normed linear space X,
subset A is bounded if and only if A is weakly bounded. Subset A
is weakly bounded if and only if ($\forall x^{*} $ $\in X^{*}$)
($\sup_{x \in A } $ $\mid x^{*} (x) \mid$ $\leq K_{x^{*}}$) where
$K_{x^{*}}$ is a constant which
depends on $x^{*}$. }} \\

{\bf Lemma 2.  Goldstine Theorem.}{ \it {For any normed linear
space X, $B_{X}$ is $weak^{*}$ dense in $B_{X^{**}}$.
Consequently, $\lambda B_{X}$ is $weak^{*}$ dense in  $\lambda
B_{X^{**}}$. (The map $T_ {\lambda}$ : x $\longrightarrow$
$\lambda x$ is a homeomorphism on convex spaces (all normed spaces
are convex spaces); since the property of being dense is a
topological
property, it is preserved under a homeomorphism.)}} \\

{\indent \bf Lemma 3.\footnote {M.G Murdeschwar {\it {General
Topology}} p72 Theorem 4.4}}{\it { Let S be a subspace of X. Let A
be any
subset of S. Then, $\bar {A}^{S} = \bar {A}^{X} \cap S $ }}\\

{\bf {Lemma 4.\footnote {D.S Bridges {\it {Foundations of Real and
Abstract Analysis}} p264 Proposition 6.1.7}}} {\it { If $x_{0}$ is
a nonzero element of a normed space X, there exists a bounded
linear functional f on X such that f($x_{0}$) = $\parallel x_{0}
\parallel$  and $\parallel f
\parallel$ = 1}}\\

{\bf {Lemma 5.}} {\it {If A is a relatively weak compact subset of
a normed space X such that $X^{*}$ contains a countable total set,
then $\bar{A}^{weak}$} is metrizable. Further the metric
restricted to A $\times$ A is a metric that generates the weak
topology of A.}

{\bf {Proof.}} Let $\Gamma $ denote the total (point separating)
subset of $X^{*}$. Let $\mathcal {J}^{X^{*}}$ and $\mathcal
{J}^{\Gamma}$ denote the subspace topology of $\bar{A}^{weak}$
(which is $\sigma (X, X^{*})$ and $\sigma (X, \Gamma)$ compact, by
hypothesis and the fact that compactness is contractive,
respectively) with respect to $\sigma (X, X^{*})$ and $\sigma (X,
\Gamma)$, respectively. Then $\mathcal {J}^{X^{*}} = \mathcal
{J}^{\Gamma}$since both are compact2 topologies. Hence, it
suffices to work in (X, $\sigma (X, \Gamma)$) to determine the
subspace topology on $\bar{A}^{weak}$.

Since $\Gamma$ is countable, (X, $\sigma (X, \Gamma)$) is
metrizable so that the subspace topology on $\bar{A}^{weak}$ is
the restriction of this metric on $\bar{A}^{weak} \times
\bar{A}^{weak}$. The identity map from $(\bar{A}^{weak}, metric \
topology)$ into $(\bar{A}^{weak}, \sigma (X, X^{*}))$ is a
continuous map from a compact space into a Hausdorff space. Hence
a homeomorphism and so A with the topology relativized with
respect to $\sigma (X, X^{*})$ is generated by the restriction of
the metric on A $\times$ A.

\begin{flushright} Q.E.D \end{flushright}
\pagebreak

 {\bf {Lemma 6.}} {\it {The dual of a separable normed
space contains a countable total set. }}

{\bf {Proof.}} Suppose X is a separable normed space. \\ \\
\indent {\bf {Step 1:}}{\it {Construct a countable dense set in
the unit sphere $B_{X}$. }}

Let $C = \{ x_{1}, x_{2}, \ldots , x_{n}, \dots \}$ denote the
countable dense subset of X. \\ Then
 \begin{equation} 
 (\forall y \in B_{X}) (\forall m \in {\bf {N}} ) (\exists x^{m} _{n} \in C
 \mid \ d(y, x^{m}_{n})) \leq 1 \backslash m)
\end{equation}

Equivalently, we have a countable index set $\Lambda$, defined as
follows (recall that y resides in both $B_{X}$ and $B_{1
\backslash m } (x^{m}_{n})$ ):
\begin{equation} 
\Lambda = \{ (n^{m}, m) \mid B_{X} \cap B_{1 \backslash m} (x^{m}
_{n}) \not= \emptyset \}
\end{equation}

For each $(n^{m}_{m}, m)$ in $\Lambda$, choosing an element,
$y^{n^{m}}_{m}$, in the non empty closed set $ B_{X} \cap B_{1
\backslash m} (x^{m} _{n}) $, we have a countable subset of
$B_{X}$ which we claim is dense in $B_{X}$. Denote this set as D.

Consider an arbitrary y $\in B_{X}$. Then for arbitrary $\epsilon
\ge 0 \ \ \exists m \in {\bf {N}} \mid \ \ 1 \backslash m \leq
\epsilon \backslash 2$.  Invoking (7), $\exists x^{m} _{n} \in C
\mid  y \in B_{1 \backslash m }(x^{m}_{n})$. But $y^{n^{m}}_{m}$
in D, is also in $B_{1 \backslash m} (x^{m} _{n})$. Hence
d$(y,y^{n^{m}}_{m})$ can be at most diameter of $B_{1 \backslash
m} (x^{m} _{n})$ apart. That is, $2 \backslash m$ or $\epsilon$ -
apart.\\

{\bf {Step 2:}}{\it{Construct a countable subset of dual $X^{*},
\{ d^{*}_{n} \}$ such that  $d^{*}_{n} (d_{n})$ = 1 where $d_{n}
\in D$ }}\\

By Lemma 4, $\forall d_{n} \in D, \exists f \in X^{*} \mid
f(d_{n}) = \parallel d_{n} \parallel$. Consequently, setting
$d^{*}_{n} = \frac{1}{\parallel d_{n} \parallel} f$ completes step
2.\\

\noindent We now show that D is a total set:-\\

\noindent ($\forall z \neq 0) \exists \lambda \ge 0 \mid \lambda z
\in B_{X}$\\ Choose $y^{n^{m}}_{m} \in D \mid$ d$(y^{n^{m}}_{m} ,
\lambda z) \le \frac{1}{2 \parallel (y^{n^{m}}_{m})^{*}
\parallel}$. \\ Now, \begin{eqnarray} \mid (y^{n^{m}}_{m})^{*}
(y^{n^{m}}_{m} - \lambda z )\mid & = & \mid 1 -
(y^{n^{m}}_{m})^{*}(\lambda z) \mid \nonumber\\
& \leq  & \parallel (y^{n^{m}}_{m})^{*} \parallel \parallel
(y^{n^{m}}_{m} - \lambda z ) \parallel \nonumber\\
& \le & \parallel (y^{n^{m}}_{m})^{*} \parallel \cdot \frac{1}{2
\parallel (y^{n^{m}}_{m})^{*} \parallel} \nonumber\\
& = & \frac{1}{2}.\nonumber
\end{eqnarray}

If $\mid 1 - (y^{n^{m}}_{m})^{*}(\lambda z) \mid  \le \frac{1}{2}$
then $(y^{n^{m}}_{m})^{*}(\lambda z) \ne 0.$ Therefore,
$(y^{n^{m}}_{m})^{*}( z) =
\frac{1}{\lambda}(y^{n^{m}}_{m})^{*}(\lambda z) \ne 0.$
\begin{flushright} Q.E.D \end{flushright}

\pagebreak
{\bf Lemma 7:} {\it {Let $E$ be a finite dimensional subspace of
the bidual $X^{**}$. Then there exists a ${\bf{finite}}$ subset $E
^{\prime}$ of the unit {\bf {sphere}} of the dual $X^{*},
S_{X^{*}}$ such that:

\begin{equation}
\frac{\parallel x^{**} \parallel}{2} \ \leq max \ \{ |
x^{**}(x^{*})| \  where \  x^{*} \in E^{\prime}  \}
\end{equation}

\noindent for each $x^{**} \in \ finite \ dimensional \ subspace \
E$. } }\\

\noindent [NOTE:$S_{X^{*}} = \{ x^{*} \in X^{*} | \parallel x^{*}
\parallel = 1 \}$ ] \\

{\bf {Lemma 8:\footnote {H.H Schaefer {\it {Topological Vector
Spaces}} 2nd Edition p26 Chapter 1 Section 5 Corollary 2}}} {\it {
Let L be a t.v.s. over K. Then the range of every
Cauchy sequence is bounded}}\\

{\bf {Lemma 9:\footnote {A. P. Robertson {\it {Topological Vector
Spaces}} p67 Chapter 4 Section 4 Theorem 1}}} {\it {The same sets
are bounded in every topology of a dual pair}}

\begin{center} ACKNOWLEDGEMENTS \end{center}

\noindent This paper would not have been written up without
Professor Elemer Rossinger. His advise on the presentation of the
paper was indispensable. I also thank Professor Anton Stroh for
the course which included the Eberlein Smulian Theorem.

\bigskip
\begin{center} REFERENCES \end{center}

\noindent 1. J. Diestel, {\it{Sequences and series in Banach
spaces}},Graduate Texts in Mathematics 92 (Springer-Verlag New
York, Inc. 1984)\\ \\

\noindent 2. M.G Murdeschwar, {\it {General Topology}}, 2nd
Edition (Wiley Eastern Limited 1990) \\ \\

\noindent 3. D.S Bridges, {\it {Foundations of Real and Abstract
Analysis}}, Graduate Texts in Mathematics 174 (Springer-Verlag New
York, Inc. 1998) \\ \\

\noindent 4. H.H Schaefer, {\it {Topological Vector Spaces}}, 2nd
Edition, Graduate Texts in Mathematics 3 (Springer-Verlag New
York, Inc. 1999) \\ \\

\noindent 5. A. P. Robertson and W. Robertson, {\it {Topological
Vector Spaces}}, Cambridge Tracts in Mathematics and Mathematical
Physics (Cambridge at the University Press 1964)\\ \\

\end{document}